%% file: main.tex
\documentclass[letterpaper, 10 pt, conference]{ieeeconf}

\IEEEoverridecommandlockouts
\AtBeginDocument{\raggedbottom}
\AtBeginDocument{\addtolength{\topmargin}{4pt}}

\usepackage{amsmath,amssymb}
\usepackage{bm}
\usepackage{booktabs}
\usepackage{multirow}
\usepackage{threeparttable}

\usepackage{graphicx}
\usepackage{cite}
\usepackage{algorithm}
\usepackage{algpseudocode}
\usepackage{tikz}
\usetikzlibrary{arrows.meta,shapes.geometric,fit,positioning,backgrounds}
\graphicspath{{figures/}}
\makeatletter\let\acc@normalsize\normalsize
\renewcommand{\normalsize}{\acc@normalsize
  \setlength{\abovedisplayskip}{4pt plus 1pt minus 1pt}
  \setlength{\belowdisplayskip}{4pt plus 1pt minus 1pt}
  \setlength{\abovedisplayshortskip}{1pt plus 1pt}
  \setlength{\belowdisplayshortskip}{2pt plus 1pt}}
\makeatother\normalsize

\newcommand{\Nhat}{\widehat N}
\newcommand{\LD}{\mathcal L_D}
\newcommand{\LM}{\mathcal L_M}
\newcommand{\Lf}{\mathcal L_f}
\newcommand{\Sset}{\mathcal S}

\title{\LARGE \bf
Learning Behavioral Parameters of Network Equilibrium Models from Aggregate Mobility Data
}

\author{Tenglong Hong$^{1}$ and Youngseo Kim$^{1}$
\thanks{$^{1}$The authors are with the Department of Civil and Environmental
        Engineering, University of California, Los Angeles, CA 90095, USA
        {\tt\small tlhong0802@g.ucla.edu}
        }
}

\begin{document}

\maketitle
\thispagestyle{empty}
\pagestyle{empty}

\begin{abstract}
Estimating travel behavior parameters is essential for understanding current travel patterns and for forecasting how travelers may respond to future changes in transportation systems. However, reliable estimation remains difficult when behavior must be inferred from sparse and aggregate observations, especially when jointly recovering several behavioral components. In this paper, we develop a bi-level estimation program for a combined destination, mode and route choice model with congestion. The inner level is a convex program at given scale parameters, and returns the taste coefficients as its dual variables together with the equilibrium travel pattern. The outer level selects the scale parameters by out-of-sample prediction, measured by how well the modeled travel pattern predicts observations withheld from estimation. To solve the bi-level problem at scale, we propose a matrix-free Hierarchical Lagrangian Dual algorithm as the inner solver and an outer search by implicit differentiation, which uses the sensitivities of the inner solution to enable efficient search over the scale parameters and reaches the calibrated scales in a handful of iterations. We validate the approach on Los Angeles County using public origin--destination, mode-share and traffic-count data. Our approach scales to a county-level problem with 2.7 million variables and substantially improves computational efficiency over a general-purpose conic solver. Compared with the fixed-scale benchmark, the calibrated model improves predictive performance at the destination and mode levels, with consistent gains across complementary goodness-of-fit measures.
\end{abstract}

\section{INTRODUCTION}
Estimating travel behavior parameters is essential for understanding current travel patterns and for forecasting how travelers may respond to future changes in transportation systems \cite{benakivaLerman1985}. However, reliable estimation remains challenging when inferring behavior at population scale. Household travel surveys provide rich individual-level information, but their inherently limited sample sizes introduce sampling uncertainty and make reliable population-wide expansion difficult. Meanwhile, increasingly available aggregate data, including origin--destination (OD) matrices, mode shares, and traffic counts, offer complementary population-level observations. These sources, however, are aggregates over different levels of the same travel decision (e.g., destination, mode, and route), so integrating them into a unified behavioral model remains challenging \cite{cascettaRusso1997}.

A recent line of work addresses this challenge by treating the network equilibrium model itself as the object to be learned: a forward problem describes the equilibrium travel pattern, and an inverse problem selects the parameters under which that pattern best reproduces the observations \cite{liuYin2026}. Within this line, the combined convex framework developed in \cite{kim2026estimate} represents destination, mode and route choices jointly with congestion, so that the equilibrium pattern and the behavioral parameters that generate it are obtained from a single program. Three challenges stand between this framework and practice, and they coincide with questions the inverse-learning perspective \cite{liuYin2026} leaves open. First, the scale parameters of the choice hierarchy are not identified by aggregate observations and remain to be calibrated. Second, the inverse problem must be solved at full real-world scale, where the forward problem alone is a county-sized equilibrium with millions of unknowns, making general-purpose solvers impractical in both memory use and computation time. Third, whether the framework can recover meaningful behavioral parameters from limited public data and deliver out-of-sample predictions remains an open empirical question. This paper addresses the three challenges, in that order. Our contributions are as follows.
\begin{itemize}
\item \emph{A bi-level estimation program.} We formulate the estimation problem as a bi-level program. The outer level is the calibration, which selects the scale parameters by the model's out-of-sample predictive accuracy. The inner level is the estimation engine: at given scale parameters, a convex program returns the equilibrium travel pattern together with the taste coefficients.
\item \emph{Scalable algorithms for both levels.} We propose novel algorithms that solve this bi-level program at county scale. The inner problem is solved by a matrix-free Hierarchical Lagrangian Dual (HLD) algorithm that scales linearly with the problem size. The outer problem is solved by a search guided by sensitivities of the inner solution, obtained through implicit differentiation.
\item \emph{A county-scale validation on public data.} We validate the approach on commuting in Los Angeles County using the LODES commute matrix \cite{lehdLodes}, ACS mode shares \cite{acs2019}, and PeMS traffic counts \cite{pems}. At a scale of 2,250 home tracts and 1.2 million routes on the road network, we evaluate the accuracy, identifiability, and predictive performance of the estimates.
\end{itemize}

\textit{Related work.} Combined destination, mode and route choice models with congestion are classical in transportation planning and go back to the entropy and equilibrium formulations \cite{wilson1967,beckmann1956}. Their parameters are conventionally estimated from individual survey records \cite{benakivaLerman1985}. Estimation from aggregate data usually draws on observations at a single level of the travel decision, link counts for the OD matrix and a route-choice coefficient \cite{yangEtAl2001}, mode choice from market shares \cite{ren2025}, and link counts and travel times for route preferences \cite{guarda2024}. Joint estimation of all three levels from aggregate data is less studied, limited to maximum likelihood on a small network \cite{abrahamsson1999}, and more recently, convex estimation on a larger network but with fixed scales \cite{kim2026estimate}. The inverse-learning perspective of \cite{liuYin2026} unifies these estimation problems as inverse problems over a forward equilibrium model, solved at the scale of a city network by a multiconvex reformulation of user equilibrium \cite{liu2026large}. In the control of transportation networks, the same equilibrium models serve as plant models for routing control \cite{kricheneEtAl2014}, and their calibration from sparse data is addressed as the inference of link latency functions from observed equilibria \cite{thaiEtAl2015}. Across the above literature, the scales are normalized or fixed by convention except in \cite{abrahamsson1999,yangEtAl2001}, the observations rarely enter at more than one level of the travel decision, and no instance reaches a metropolitan region on its full road network. Our paper calibrates the scales, integrates observations from all three levels, and solves the problem for Los Angeles County, a metropolitan region.
\section{MODEL}
\label{sec:model}

\textit{Notation.} We use calligraphic symbols to denote sets, bold symbols to denote vectors, and a hat ( $\widehat{\cdot}$ ) to mark an observed quantity or a given parameter. Let $i\in\mathcal I$ denote origins, partitioned into estimation origins $\mathcal I_{\mathrm{est}}$, on which the model is fitted, and validation origins $\mathcal I_{\mathrm{val}}$, on which its losses are evaluated. From origin $i$, a trip chooses destination $j\in\mathcal J_i$, then mode $m\in\mathcal M_{ij}$, then route $r\in\mathcal R_{ijm}$ on a network with links $a\in\mathcal A$, where $\widehat\delta^a_{ijmr}=1$ if route $r$ uses link $a$, and each link carries flows $f^m_a$ per mode. The data are trip counts: $\widehat T_{ij}$ trips from $i$ to $j$ and $\widehat T_{ijm}$ their split by mode $m$, with productions $\widehat O_i=\sum_j\widehat T_{ij}$ and total $\Nhat=\sum_i\widehat O_i$. Divided by $\Nhat$, they give the empirical probabilities $\widehat p_i$, $\widehat p_{ij}$ and $\widehat p_{ijm}$. The data also include the attribute vectors $\widehat{\mathbf X}^D_{ij}$ for an OD pair and $\widehat{\mathbf X}^M_{ijm}$ for mode $m$ on that pair.

The travel pattern is described by the joint probabilities $p_{ij}$, $p_{ijm}$ and $p_{ijmr}$ that a trip has a given destination, mode and route, equivalently the trip flows $T_{ij}=\Nhat p_{ij}$, $T_{ijm}=\Nhat p_{ijm}$ and $T_{ijmr}=\Nhat p_{ijmr}$, which in turn induce the link flows $f^m_a$ through $f^m_a=\sum_{i,j,r}T_{ijmr}\widehat\delta^a_{ijmr}$. The ratios of the probabilities are the conditional choice probabilities $P_D(j\mid i)=p_{ij}/\widehat p_i$, $P_M(m\mid i,j)=p_{ijm}/p_{ij}$ and $P_R(r\mid i,j,m)=p_{ijmr}/p_{ijm}$. The taste coefficients are $\boldsymbol\beta^D$ and $\boldsymbol\beta^M$, associated with the attributes of the destination and mode levels. The scale parameters are $\mathbf s=(\theta_D,\theta_M,\lambda)$, where $\theta_D$ and $\theta_M$ set how sharply the destination and mode levels respond to utility differences, while $\lambda$ governs how strongly route choice responds to congested travel time.

\subsection{The choice hierarchy}
\label{subsec:hierarchy}

We represent travel demand as a hierarchy of choices in equilibrium with the network. Each traveler departing from origin $i$ chooses destination $j\in\mathcal J_i$, mode $m\in\mathcal M_{ij}$ and route $r\in\mathcal R_{ijm}$, and starting from the $\widehat O_i$ trips produced at $i$, the three conditional probabilities allocate trips level by level into mode-specific link flows,
\begin{equation*}
    \widehat O_i
    \xrightarrow{\;P_D\;} T_{ij}
    \xrightarrow{\;P_M\;} T_{ijm}
    \xrightarrow{\;P_R\;} T_{ijmr}
    \xrightarrow{\;\widehat\delta^a_{ijmr}\;} f_a^m .
\end{equation*}
The link flows $\mathbf f=(f_a^m)_{a,m}$ set the congested travel times, which enter the utilities of all three levels, so a fixed point of this mapping is an equilibrium in which travelers respond to the congestion their own choices generate. Fig.~\ref{fig:choice-hierarchy} previews the hierarchy and the bi-level estimation framework, whose elements are defined in the rest of this section.

\input{inputs/bilevel_framework}

\textit{Destination utility.} Destination choice follows multinomial logit \cite{benakivaLerman1985}, which is equivalent to a gravity model \cite{anas1983}, with utility $V^D_{ij}=\boldsymbol\beta^{D\top}\widehat{\mathbf X}^D_{ij}$, where $\widehat{\mathbf X}^D_{ij}=\big(\ln(1+\widehat A_j),\ \widehat d_{ij},\ \ln\widehat d_{ij},\ \widehat\kappa_{ij}\big),$ with corresponding coefficients $
\boldsymbol\beta^D=\big(\beta^D_{\mathrm{size}},\beta^D_{\mathrm{dist}},\beta^D_{\log},\beta^D_{\mathrm{sec}}\big).$
Here $\widehat A_j$ is the number of jobs at destination $j$, which often enters in logarithm. The distance $\widehat d_{ij}$ enters through the distance-deterrence function with two parameters $\gamma$ and $\eta$, $d^{\gamma}e^{\eta d}, $ as commonly used in gravity literature \cite{lenormandEtAl2016}. Its logarithm yields the two terms $\beta^D_{\log}\ln\widehat d_{ij}+\beta^D_{\mathrm{dist}}\widehat d_{ij}$ with $\beta^D_{\log}=\gamma$ and $\beta^D_{\mathrm{dist}}=\eta$. The sector match $\widehat\kappa_{ij}=\sum_s\widehat w_{is}\widehat q_{js}$ is the probability that a job at $j$ is in the sector of a worker living at $i$, computed from the share $\widehat w_{is}$ of $i$'s resident workers and the share $\widehat q_{js}$ of $j$'s jobs in each industry sector $s$.

\textit{Mode utility.} Mode choice follows a nested logit \cite{benakivaLerman1985}. The mode set $\mathcal M_{ij}$ is partitioned into nests $\mathcal M_{ijn}$, $n\in\mathcal N$, and a mode $m\in\mathcal M_{ijn}$ has utility $V^M_{ijm}=\boldsymbol\beta^{M\top}\widehat{\mathbf X}^M_{ijm}$, where $\widehat{\mathbf X}^M_{ijm}=\big(\mathbf 1_m,\ \widehat t_{ijm}\big)$ and $\boldsymbol\beta^M=\big(\boldsymbol{\mathrm{ASC}},\ \beta^M_{\mathrm{time}}\big)$. Here $\mathbf 1_m$ denotes the indicator for mode $m$, with driving alone as the reference alternative whose alternative-specific constant is normalized to zero, and $\widehat t_{ijm}$ is the door-to-door travel time, so that the utility reads $\mathrm{ASC}_m+\beta^M_{\mathrm{time}}\widehat t_{ijm}$.

\textit{Route utility.} Route choice follows path-size logit \cite{benakivaBierlaire1999} with utility $V^R_{ijmr}=-\lambda\,g_{ijmr}(\mathbf f)+\ln\widehat\psi_{ijmr}$. The congested route time is $g_{ijmr}(\mathbf f)=\sum_a\widehat\delta^a_{ijmr}\,g^m_a(f^m_a)$, where $g^m_a$ is the Bureau of Public Roads (BPR) latency $g^m_a =  \widehat t^{\,0}_{am}\big[1+\widehat a^m(f^m_a/\widehat c^m_a)^{\widehat b^m}\big]$ with free-flow time $\widehat t^{\,0}_{am}$, capacity $\widehat c^m_a$ and shape parameters $ \widehat a^m$, $\widehat b^m$, and $\widehat\psi_{ijmr}$ is the path-size factor that corrects for overlapping routes. In any hierarchical logit, the scale of one level is fixed for identification and the others are measured relative to it \cite{benakivaLerman1985}. We fix the route level's scale at 1, as in \cite{kim2026estimate}, so that route utilities are measured in utils. The weight $\lambda$, in utils per hour, converts the congested route time $g_{ijmr}(\mathbf f)$ from hours into utils, and a larger $\lambda$ makes route choice respond more strongly to travel time.

Given the utilities, the choice probabilities follow the logit form at each level, and the levels are linked from the bottom up: a mode inherits the expected maximum utility of its routes through the logsum $S^R_{ijm}=\ln\sum_{r\in\mathcal R_{ijm}}e^{V^R_{ijmr}}$, a nest $n\in\mathcal N$ carries a dissimilarity parameter $\widehat\tau_n\in(0,1]$ and inherits that of its modes through $S^N_{ijn}=\ln\sum_{m\in\mathcal M_{ijn}}e^{\theta_M U_{ijm}/\widehat\tau_n}$ with $U_{ijm}=V^M_{ijm}+S^R_{ijm}$, and a destination inherits that of its nests through $S^M_{ij}=\theta_M^{-1}\ln\sum_{n\in\mathcal N}e^{\widehat\tau_n S^N_{ijn}}$. The conditional probabilities are: 
\begin{align}
    P_R(r\mid i,j,m)&=e^{V^R_{ijmr}}\Big/\textstyle\sum_{r'\in\mathcal R_{ijm}}e^{V^R_{ijmr'}},
    \label{eq:route}\\
    P_M(m\mid i,j)&=\frac{e^{\theta_M U_{ijm}/\widehat\tau_n}\big(\sum_{m'\in\mathcal M_{ijn}}e^{\theta_M U_{ijm'}/\widehat\tau_n}\big)^{\widehat\tau_n-1}}{\sum_{n'\in\mathcal N}\big(\sum_{m'\in\mathcal M_{ijn'}}e^{\theta_M U_{ijm'}/\widehat\tau_{n'}}\big)^{\widehat\tau_{n'}}}
    \label{eq:mode}\\
    P_D(j\mid i)= &e^{\theta_D(V^D_{ij}+S^M_{ij})}\Big/\textstyle\sum_{j'\in\mathcal J_i}e^{\theta_D(V^D_{ij'}+S^M_{ij'})},
    \label{eq:destination}
\end{align}
with $n$ the nest containing $m$. Two properties of the scales matter here. First, only the products $\theta_D\boldsymbol\beta^D$ and $\theta_M\boldsymbol\beta^M$ enter \eqref{eq:mode}--\eqref{eq:destination}, so a taste coefficient is comparable across scales only through its product. Second, consistency of the hierarchy requires $1\ge\theta_M\ge\theta_D$. A scale is inversely related to the variance of the unobserved utility at its level, and that variance accumulates going up the hierarchy, so the scale cannot increase \cite{benakivaLerman1985}.

\subsection{Inner problem: estimation as maximum entropy}
\label{subsec:inner}

The inner problem is the estimation engine. At given scales $\mathbf s$, it returns the equilibrium travel pattern together with the taste coefficients $\boldsymbol\beta^D$ and $\boldsymbol\beta^M$, using the data of the estimation origins $\mathcal I_{\mathrm{est}}$ only. Its variables are the travel pattern, the joint probabilities $\mathbf p=(p_{ij},p_{ijm},p_{ijmr})$ and the link flows $\mathbf f$, constrained to the sets
\begin{subequations}
\begin{align}
    \mathcal P&=\big\{\mathbf p\ge\mathbf 0:\ \textstyle\sum_j p_{ij}=\widehat p_i,\ \sum_m p_{ijm}=p_{ij}, \notag \\
    &\hphantom{=\big\{\mathbf p\ge\mathbf 0:\ }\textstyle\sum_r p_{ijmr}=p_{ijm}\big\},\label{eq:pattern-set}\\
    \mathcal F(\mathbf p)&=\big\{\mathbf f:\ f^m_a/\Nhat=\textstyle\sum_{i,j,r}p_{ijmr}\widehat\delta^a_{ijmr}\ \ \forall\,a,m\big\},\label{eq:route-link}
\end{align}
\label{eq:consistency}
\end{subequations}
where $\mathcal P$ collects the hierarchical consistency and adding-up conditions and $\mathcal F(\mathbf p)$ the route--link condition that loads the route pattern onto the network links.

The observations enter through moment constraints. For each attribute of the destination and mode levels, the modeled trips must reproduce the observed average of that attribute. These are the natural estimating conditions for a logit model, since the first-order conditions of maximum likelihood equate the model and observed averages exactly \cite{benakivaLerman1985}. Most usefully, the multiplier of each constraint is the taste coefficient of its attribute, $\boldsymbol\beta^D$ for the destination moments and $\boldsymbol\beta^M$ for the mode moments, which links the moment-matching primal to the behavioral parameters. The two moment constraints are therefore
\begin{equation}
\begin{aligned}
  \textstyle\sum_{i,j}p_{ij}\,\widehat{\mathbf X}^D_{ij}&=\textstyle\sum_{i,j}\widehat p_{ij}\,\widehat{\mathbf X}^D_{ij}
    &&[\boldsymbol\beta^D],\\
    \textstyle\sum_{i,j,m}p_{ijm}\,\widehat{\mathbf X}^M_{ijm}&=\textstyle\sum_{i,j,m}\widehat p_{ijm}\,\widehat{\mathbf X}^M_{ijm}
    &&[\boldsymbol\beta^M].
\end{aligned}
    \label{eq:moments}
\end{equation}
with the multiplier of each block in brackets, $\boldsymbol\beta^D$ for the destination moments and $\boldsymbol\beta^M$ for the mode moments. The moment formulation allows aggregate observations at different levels of the travel hierarchy to be integrated through their corresponding moments, while the available moments determine which parameters can be recovered. 

The objective follows the maximum-entropy principle: among all travel patterns that satisfy the hierarchical consistency and moment constraints, the inner problem selects the one that assumes the least beyond the data, which is the pattern of largest entropy \cite{wilson1967}. We define the objective function of the inner problem as
\begin{equation}
\begin{aligned}
    J(\mathbf p,\mathbf f)={}&\tfrac{1}{\theta_D}H_D(\mathbf p)+\tfrac{1}{\theta_M}H_M(\mathbf p)+H_R(\mathbf p)\\
    &-\tfrac{\lambda}{\Nhat}\textstyle\sum_{a,m}\int_0^{f_a^m}g_a^m(w)\,dw,
\end{aligned}
    \label{eq:objective}
\end{equation}
with the conditional entropies $H_D=-\sum_{ij}p_{ij}\ln(p_{ij}/\widehat p_i)$, $H_M=-\sum_{ijm}p_{ijm}\ln(p_{ijm}/p_{ij})$ and $H_R=-\sum_{ijmr}p_{ijmr}\ln\big(p_{ijmr}/(\widehat\psi_{ijmr}p_{ijm})\big)$. For our hierarchical model, the total entropy is a sum of one conditional entropy per level, each weighted by the reciprocal of the level's scale, so a smaller scale gives its level more weight in the objective and leaves it more dispersed. The last term is the Beckmann potential of the network \cite{beckmann1956}. The network potential sums over links the latency function integrated from $0$ to each link's flow $f_a^m$, so its derivative with respect to a link flow is the corresponding congested travel time. Because the derivative of the potential is the congested time, the optimality condition in $p_{ijmr}$ reproduces the route logit \eqref{eq:route} at the congested times of the optimal flows. The optimal flows are therefore consistent with the route choice they induce, which is the definition of a stochastic user equilibrium. Collecting these elements, the inner problem given $\mathbf s$ is
\begin{equation}
\label{eq:inner}
\boxed{\;
\begin{aligned}
    &\textsc{Inner Problem}\ \big(\text{given} \ \mathbf s\big)\\
    \max_{\mathbf p\ge\mathbf 0,\,\mathbf f}\quad & J(\mathbf p,\mathbf f)\ \text{in \eqref{eq:objective}}\\
    \text{s.t.}\quad & \text{consistency conditions \eqref{eq:consistency}},\\
    & \text{moment constraints \eqref{eq:moments}}\quad[\boldsymbol\beta^D,\boldsymbol\beta^M].
\end{aligned}\;}
\end{equation}
Under the standard assumptions that the latencies $g^m_a$ are increasing in the flow and the scales are positive, \eqref{eq:inner} is a convex program, and its optimality conditions are exactly \eqref{eq:route}--\eqref{eq:destination} with the taste coefficients recovered as the multipliers of the moment constraints \eqref{eq:moments} \cite{kim2026estimate}. It thus extends the logit estimation principle to the full choice hierarchy: sensitivities ordinarily estimated from individual survey records are recovered from public aggregate totals alone, and the joint estimation across destination, mode, route and congestion lets each level inform the others, for internally consistent estimates. Solving the inner problem at given scales therefore returns the equilibrium pattern and the taste coefficients at once.

\subsection{Outer problem: scales selected by prediction}
\label{subsec:outer}

The outer level selects the scales by out-of-sample prediction. For a candidate $\mathbf s=(\theta_D,\theta_M,\lambda)$, the inner problem is solved on $\mathcal I_{\mathrm{est}}$ and yields the taste coefficients $\boldsymbol\beta^\star(\mathbf s)$. With these coefficients fixed, the model predicts the travel pattern of the validation origins $\mathcal I_{\mathrm{val}}$ from their observed productions and choice sets. The predictions are compared with three observations the inner problem does not use, the destination flows and mode shares of the validation origins and the traffic counts on the network, through the following three losses.

\textit{Destination loss.} Let $\widehat P_D(j\mid i)=\widehat T_{ij}/\widehat O_i$ be the observed share of commuters who live in origin $i$ and work in destination $j$, and $\widehat{w}_i=\widehat O_i/\sum_{i'\in\mathcal I_{\mathrm{val}}}\widehat O_{i'}$ the share of commuters from the validation origin $i$. The destination loss is the cross-entropy between the observed and the predicted shares, averaged over the validation origins with $\widehat{w}_i$,
\begin{equation*}
    \LD=-\sum_{i\in\mathcal I_{\mathrm{val}}}\widehat{w}_i\sum_{j\in\mathcal J_i}\widehat P_D(j\mid i)\,\ln P_D(j\mid i;\mathbf s,\boldsymbol\beta^\star).
\end{equation*}
The log loss of a predicted distribution against an observed one is the entropy of the observed distribution plus the Kullback--Leibler (KL) divergence between the two, so $\LD=\widehat{\mathcal H}_{\mathrm{val}}+\mathrm{KL}_{\mathrm{val}}$, with $\mathrm{KL}_{\mathrm{val}}=\sum_{i\in\mathcal I_{\mathrm{val}}}\widehat{w}_i\sum_{j\in\mathcal J_i}\widehat P_D(j\mid i)\ln\frac{\widehat P_D(j\mid i)}{P_D(j\mid i)}\ \ge0,$
and $\widehat{\mathcal H}_{\mathrm{val}}=-\sum_i \widehat{w}_i\sum_j\widehat P_D(j\mid i)\ln\widehat P_D(j\mid i)$ which depends on the data alone. $\widehat{\mathcal H}_{\mathrm{val}}$ is the minimum attainable log loss on this sample, achieved when the predicted shares exactly match the observed shares. $\mathrm{KL}_{\mathrm{val}}$ measures the additional loss due to discrepancies between the predicted and observed shares.

\textit{Mode loss.} The mode loss compares the model-implied mode split of a validation origin $i$, $\sigma_{im}=\sum_j P_D(j\mid i)\,P_M(m\mid i,j)$, with the observed split $\widehat\sigma_{im}=\sum_j\widehat T_{ijm}/\widehat O_i$, for each validation origin $i$. Mode shares are conventionally stated in percent, so the loss is expressed in percentage points, the mean absolute difference of the shares multiplied by 100, $\LM=\frac{100}{|\mathcal I_{\mathrm{val}}|\,|\mathcal M|}\sum_{i\in\mathcal I_{\mathrm{val}}}\sum_{m\in\mathcal M}\big|\sigma_{im}-\widehat\sigma_{im}\big|.$

\textit{Flow loss.} The flow loss compares the model with the traffic count $\widehat f_d$ at each detector $d\in\mathcal D$. A detector observes a road segment that the network represents as a set of consecutive links, $\mathcal A_d$, so the modeled flow at the detector is obtained by summing over the road-using modes and averaging over those links, $f_d=\frac{1}{|\mathcal A_d|}\sum_{a\in\mathcal A_d}\sum_{m}f^m_a$. The traffic counts and the modeled flow can differ by a scale factor $\nu$, since the estimation data may cover only a subset of trip purposes and thus not represent the full traffic volume counted by the detectors. In logarithms, a common factor becomes a common shift, so the log ratio $\rho_d=\ln f_d-\ln\widehat f_d$ carries the same shift $\ln\nu$ at every detector. By subtracting the mean of the log ratios $\bar\rho=\frac{1}{|\mathcal D|}\sum_{d\in\mathcal D}\rho_d$ over detectors, we remove the common scale difference $\ln\nu$ so the loss compares the relative distribution of the flows across detectors and is independent of their overall level. The loss is therefore the root-mean-square error of the log flows about that mean, $\Lf=\Big(\frac{1}{|\mathcal D|}\sum_{d\in\mathcal D}(\rho_d-\bar\rho)^2\Big)^{1/2}.$

The three losses are in different units, so we combine them by the $\varepsilon$-constraint method of multi-objective optimization \cite{haimes1971}. The destination loss is minimized, because the origin--destination flows are the most reliable observation in our data, and the mode and flow losses are bounded by the tolerances $\varepsilon_M$ and $\varepsilon_f$ in their own units. The roles can be reassigned when another observation is the most reliable. The scales range over $\Sset=\{\mathbf s:\ 0<\theta_D\le\theta_M\le1,\ \lambda>0\}$, as required by hierarchical consistency. The bi-level program is
\begin{equation}
\label{eq:bilevel}
\boxed{\;
\begin{aligned}
    &\textsc{Bi-level Estimation Problem}\\
    \min_{\mathbf s\in\Sset}\quad & \LD\big(\boldsymbol\beta^\star(\mathbf s)\big)\\
    \text{s.t.}\quad & \LM\big(\boldsymbol\beta^\star(\mathbf s)\big)\le\varepsilon_M,\quad \Lf\big(\boldsymbol\beta^\star(\mathbf s)\big)\le\varepsilon_f,\\
    & (\mathbf p^\star,\mathbf f^\star)\ \text{solves the inner problem \eqref{eq:inner} at }\mathbf s,\\
    & \boldsymbol\beta^\star(\mathbf s)\ \text{its multipliers on \eqref{eq:moments}}.
\end{aligned}\;}
\end{equation}

Fig.~\ref{fig:choice-hierarchy} summarizes the resulting bi-level structure.

\section{ALGORITHMS}
\label{sec:algorithm}

Bi-level programs are computationally challenging in general, since the outer objective depends on the inner solution and can be nonconvex in the outer variables even when the inner problem is convex. Our program has a structure that we leverage: the outer level calibrates only three scale parameters, whereas the inner level carries millions of variables at county scale and must be solved at every outer iterate. We therefore design an efficient inner solver first, as a fast and exact oracle, and build the outer search on it.

\subsection{Inner solver: Hierarchical Lagrangian Dual (HLD)}
\label{subsec:hld}

The constraints of \eqref{eq:inner} are of two kinds: (i) the conditions defining the pattern set $\mathcal P$ are separable across origins, each involving the probabilities of one origin alone; and (ii) the moment constraints \eqref{eq:moments} couple the origins to each other, and the route--link condition \eqref{eq:route-link} couples them to the flows as well. We dualize only these coupling constraints, so that the choice problem decomposes origin by origin in the primal. The variables of the resulting partial dual are only the taste coefficients and one multiplier per mode-specific link, far fewer than the millions of pattern variables, and its gradient is the vector of residuals of the dualized constraints at the primal maximizer. Since the primal maximizer is already obtained in closed form when evaluating the dual objective, these residuals can be computed for free. We show the details in \eqref{eq:gradient} below.

\textit{Dual function.} Let $\boldsymbol\beta^D$ and $\boldsymbol\beta^M$ be the multipliers of the moment constraints \eqref{eq:moments}, and let $\varpi^m_a$ be the multiplier of the route--link condition \eqref{eq:route-link}, one per link $a$ and mode $m$. We collect the dual multipliers in $\mathbf z=(\boldsymbol\beta^D,\boldsymbol\beta^M,\boldsymbol\varpi)$. Let $\mathbf m_D(\mathbf p)=\sum_{i,j}p_{ij}\widehat{\mathbf X}^D_{ij}$ and $\mathbf m_M(\mathbf p)=\sum_{i,j,m}p_{ijm}\widehat{\mathbf X}^M_{ijm}$ for the model's moments, $\widehat{\mathbf m}_D$, $\widehat{\mathbf m}_M$ for the observed ones, and $c_{ijmr}(\boldsymbol\varpi)=\sum_a\widehat\delta^a_{ijmr}\varpi^m_a$ for the cost of a route at the multipliers. The Lagrangian maximized over the primal variables gives the dual function as the sum of a choice block $\mathcal V(\mathbf z)$, a congestion block $\Phi(\mathbf z)$ and a data term $C(\mathbf z)$.
The resulting partial dual objective is
\begin{align}
    \mathcal D(\mathbf z)&=\mathcal V(\mathbf z)+\Phi(\mathbf z)-C(\mathbf z),\label{eq:dual}\\
    \mathcal V(\mathbf z)&=\max_{\mathbf p\in\mathcal P}\Big[\tfrac{1}{\theta_D}H_D+\tfrac{1}{\theta_M}H_M+H_R+\boldsymbol\beta^{D\top}\mathbf m_D\notag\\
    &\hphantom{=\max_{\mathbf p\in\mathcal P}\Big[}+\boldsymbol\beta^{M\top}\mathbf m_M-\textstyle\sum_{i,j,m,r}p_{ijmr}\,c_{ijmr}(\boldsymbol\varpi)\Big],\notag\\
    \Phi(\mathbf z)&=\textstyle\sum_{a,m}\ \max_{f\ge0}\Big[\dfrac{\varpi^m_a f}{\Nhat}-\dfrac{\lambda}{\Nhat}\displaystyle\int_0^{f}g^m_a(w)\,dw\Big],\notag\\
    C(\mathbf z)&=\boldsymbol\beta^{D\top}\widehat{\mathbf m}_D+\boldsymbol\beta^{M\top}\widehat{\mathbf m}_M.\notag
\end{align}
The choice block $\mathcal V(\mathbf z)$ involves only $\mathbf p$ and the congestion block $\Phi(\mathbf z)$ only $\mathbf f$, so the Lagrangian separates and each block is maximized on its own. The data term $C$ is linear in the multipliers $\boldsymbol\beta^D$ and $\boldsymbol\beta^M$.

At the maximizer of the congestion block, $\varpi^m_a=\lambda\,g^m_a(f^\star)$, so the multiplier is the congested time of the link converted into utils, and it never falls below the free-flow floor $\lambda\widehat t^{\,0}_{am}$. We therefore write $\varpi^m_a=\lambda\widehat t^{\,0}_{am}+\xi^m_a$ with an excess cost $\xi^m_a\ge0$. The excess cost is zero on links whose latency is constant in the flow, so the dual variables reduce to $\mathbf z=(\boldsymbol\beta^D,\boldsymbol\beta^M,\boldsymbol\xi)$ with $\boldsymbol\xi=(\xi^m_a)_{(a,m)\in\mathcal E_C}$ over the congestible link states $\mathcal E_C=\{(a,m):\,g^m_a\text{ is increasing}\}$, and the sum in $\Phi$ runs over $\mathcal E_C$ alone.

\textit{Dual program.} The Hierarchical Lagrangian Dual of \eqref{eq:inner} is
\begin{equation}
    \boxed{\;\min_{\mathbf z}\ \mathcal D(\mathbf z)\quad\text{s.t.}\quad \xi^m_a\ge0,\ \ (a,m)\in\mathcal E_C\;}
    \label{eq:dual-problem}
\end{equation}
Under a relative-interior constraint qualification, strong duality holds, and at a dual minimizer $\mathbf z^\star=(\boldsymbol\beta^{D\star},\boldsymbol\beta^{M\star},\boldsymbol\xi^\star)$ the block maximizers are the unique solution of \eqref{eq:inner}, with $\boldsymbol\beta^{D\star}$ and $\boldsymbol\beta^{M\star}$ its taste coefficients. The dual is easier to solve than the primal because it is convex, as a maximum of functions affine in $\mathbf z$, and differentiable. Let $\mathsf L(\mathbf p,\mathbf f;\mathbf z)$ be the Lagrangian of \eqref{eq:inner}, so that $\mathcal D(\mathbf z)=\max_{\mathbf p\in\mathcal P,\,\mathbf f\ge\mathbf 0}\mathsf L(\mathbf p,\mathbf f;\mathbf z)$ by \eqref{eq:dual}. The maximizer is unique, since the entropy is strictly concave in $\mathbf p$ and the Beckmann integral strictly convex in $\mathbf f$. By Danskin's theorem \cite{danskin1967}, the gradient of the dual is then the gradient of the Lagrangian in $\mathbf z$ at that maximizer, $\nabla_{\mathbf z}\mathcal D(\mathbf z)=\nabla_{\mathbf z}\mathsf L\big(\mathbf p^\star(\mathbf z),\mathbf f^\star(\mathbf z);\mathbf z\big)$. The partial derivative of $\mathsf L$ in each multiplier is the constraint it is attached to, so the gradient is the vector of constraint residuals at the block maximizers $\mathbf p^\star(\mathbf z)$ and $\mathbf f^\star(\mathbf z)$,
\begin{align}
    \nabla_{\boldsymbol\beta^D}\mathcal D&=\mathbf m_D(\mathbf p^\star)-\widehat{\mathbf m}_D,\quad \nabla_{\boldsymbol\beta^M}\mathcal D=\mathbf m_M(\mathbf p^\star)-\widehat{\mathbf m}_M,\notag\\
    \nabla_{\xi^m_a}\mathcal D&=\big(f^{\star m}_a(\mathbf z)-f^{\mathrm{load},m}_a(\mathbf z)\big)/\Nhat ,\label{eq:gradient}
\end{align}
where $f^{\mathrm{load},m}_a(\mathbf z)=\Nhat\sum_{i,j,r}p^\star_{ijmr}(\mathbf z)\,\widehat\delta^a_{ijmr}$ is the flow that $\mathbf p^\star(\mathbf z)$ loads on the links. A zero gradient means primal feasibility, since each block of \eqref{eq:gradient} is the residual of a dualized constraint. Additionally, since the dual is convex, it also implies dual optimality. The dual is therefore minimized by a gradient method once $\mathcal D$ and $\nabla\mathcal D$ can be evaluated.

\textit{The oracle.} What remains is to evaluate $\mathcal D$ and $\nabla\mathcal D$ at a given $\mathbf z$. Both evaluations only require the two block maximizers $\mathbf p^\star(\mathbf z)$ and $\mathbf f^\star(\mathbf z)$, which we obtain in closed form from the choice and congestion blocks without solving a new optimization problem. First, the choice block $\mathcal V(\mathbf z)$ has the same entropy-regularized hierarchical structure as the choice model of Section~\ref{sec:model}, with the multipliers $\boldsymbol\beta^D$ and $\boldsymbol\beta^M$ as the taste coefficients and the multiplier cost $c_{ijmr}(\boldsymbol\varpi)$ in place of the congested time $\lambda g_{ijmr}(\mathbf f)$. Its maximizer is therefore recovered by evaluating the nested logit \eqref{eq:route}--\eqref{eq:destination} at these utilities, as for \eqref{eq:inner} \cite{kim2026estimate}. A bottom-up pass computes the conditional probabilities and logsums, and a top-down pass chains them into $\mathbf p^\star(\mathbf z)$ and the link loads $\mathbf f^{\mathrm{load}}(\mathbf z)$.

Second, for the congestion block $\Phi(\mathbf z)$, because the Beckmann term is separable across links and modes, the maximization over the link flows decomposes by $(a,m)$, each term being the convex conjugate of the corresponding Beckmann integral. Each scalar problem is strictly concave in its link flow, since the latency is increasing, and its first-order condition equates the link's dual cost with the weighted latency at the maximizing flow, $\varpi^m_a=\lambda g^m_a(f)$. Solving this condition for $f$ under the BPR form gives $\mathbf f^\star(\mathbf z)$ explicitly,
\begin{equation}
    f^{\star m}_a(\mathbf z)=\widehat c^m_a\Big[\frac{\varpi^m_a/(\lambda\widehat t^{\,0}_{am})-1}{\widehat a^m}\Big]_+^{1/\widehat b^m},
    \label{eq:inverse-bpr}
\end{equation}
one scalar calculation per link state $(a,m)$. The oracle is therefore matrix-free.

\textit{Iteration.} With the oracle in hand, \eqref{eq:dual-problem} is solved with gradient information alone. We use L-BFGS-B \cite{byrd1995lbfgsb}, a quasi-Newton method that handles the bound $\boldsymbol\xi\ge0$ natively and, instead of forming a Hessian, approximates the curvature of $\mathcal D$ from the differences of a few recent iterates and gradients. Each quasi-Newton step is safeguarded by a descent test, and the iteration can converge globally to the dual minimizer. It stops when every residual block of \eqref{eq:gradient} is within tolerance. Algorithm~\ref{alg:hld} summarizes the method.

\begin{algorithm}[!t]
\caption{Inner solver: HLD with L-BFGS-B}
\label{alg:hld}
\begin{algorithmic}[1]
\Require scales $\mathbf s$, data of \eqref{eq:inner}, start $\mathbf z_0$, tolerance
\Repeat
    \State Route costs $c_{ijmr}(\boldsymbol\varpi)$ at $\varpi^m_a=\lambda\widehat t^{\,0}_{am}+\xi^m_a$
    \State Bottom-up pass for the logsums and the conditionals \eqref{eq:route}--\eqref{eq:destination}, top-down pass for $\mathbf p^\star(\mathbf z)$ and $\mathbf f^{\mathrm{load}}(\mathbf z)$
    \State Link flows $\mathbf f^\star(\mathbf z)$ by \eqref{eq:inverse-bpr}, then $\mathcal D$ by \eqref{eq:dual}, $\nabla\mathcal D$ by \eqref{eq:gradient}
    \State One L-BFGS-B step in $\mathbf z$ subject to $\boldsymbol\xi\ge0$
\Until{every block of \eqref{eq:gradient} is below tolerance}
\Ensure $\boldsymbol\beta^\star$, $\boldsymbol\xi^\star$, the pattern $\mathbf p^\star(\mathbf z^\star)$ and the flows $\mathbf f^\star(\mathbf z^\star)$
\end{algorithmic}
\end{algorithm}

\subsection{Outer search by implicit differentiation}
\label{subsec:outer-search}

Each loss of \eqref{eq:bilevel} depends on the scales directly and through the inner solution, $\mathcal L(\mathbf s,\mathbf z^\star(\mathbf s))$ for $\mathcal L\in\{\LD,\LM,\Lf\}$, so every evaluation costs an inner solve. Derivative-free search such as grid search or random sampling is therefore too costly at county scale, and we descend along the gradients of the losses instead. These require the derivative $\mathrm d\mathbf z^\star/\mathrm d\mathbf s$ of the inner solution, which we obtain by implicit differentiation of the optimality conditions of \eqref{eq:dual-problem}, the classical device of descent methods for bi-level programs \cite{ghadimiWang2018}. Each iteration solves a local model of \eqref{eq:bilevel} within a trust region and warm-starts the next inner solve from a first-order prediction of $\mathbf z^\star$, so an outer iteration costs a fraction of a full inner solve.

\textit{Sensitivities.} Let $\mathbf r(\mathbf z,\mathbf s)$ be the gradient \eqref{eq:gradient} restricted to $\boldsymbol\beta$ and the congested links, $\xi^m_a>0$, the coordinates of $\mathbf z$ that move with $\mathbf s$. At the minimizer of \eqref{eq:dual-problem}, $\mathbf r(\mathbf z^\star(\mathbf s),\mathbf s)=\mathbf 0$ for every $\mathbf s$. Differentiating this identity in $\mathbf s$ gives the sensitivity matrix $\mathbf U=\mathrm d\mathbf z^\star/\mathrm d\mathbf s$, whose three columns are the derivatives of $\mathbf z^\star$ with respect to $\theta_D$, $\theta_M$ and $\lambda$, as the solution of the linear system $\nabla_{\mathbf z}\mathbf r\,\mathbf U=-\nabla_{\mathbf s}\mathbf r$. We  solve this linear system by GMRES \cite{saad1986gmres}, a matrix-free approach.

\textit{Local approximation and step.} Each iteration of the outer search builds a local approximation of \eqref{eq:bilevel} around the current scales, and takes the solution of this approximate problem as the descent step for $\mathbf s$. The step $\mathbf d\in\mathbb R^3$ is taken in the logarithms of the scales, $\mathbf s\mapsto\mathbf s\,\mathrm e^{\mathbf d}$, which keeps the scales positive. Along the step, the chain rule gives the gradient of a loss at $\mathbf d=\mathbf 0$ as $\nabla_{\mathbf d}\mathcal L=\operatorname{diag}(\mathbf s)\big(\partial_{\mathbf s}\mathcal L+\mathbf U^{\top}\partial_{\mathbf z}\mathcal L\big)$, with both partial derivatives in closed form from the choice model. Applying it to $\LM$, to $\Lf$ and to $\ln P_D(j\mid i)$ for each validation pair gives $\mathbf g_M$, $\mathbf g_f$ and $\mathbf b_{ij}$. The constrained losses enter the local model linearized, $\LM+\mathbf g_M^{\top}\mathbf d$ and $\Lf+\mathbf g_f^{\top}\mathbf d$. The destination loss is linearized through its log-probabilities, $a_{ij}(\mathbf d)=\ln P_D(j\mid i)+\mathbf b_{ij}^{\top}\mathbf d$, with $\mathrm e^{a_{ij}(\mathbf d)}$ renormalized over the destinations of each origin. The local model is
\begin{equation}
\begin{aligned}
\min_{\mathbf d}\ &\sum_{i\in\mathcal I_{\mathrm{val}}}\widehat{w}_i\Big[\ln\sum_{j}\mathrm e^{a_{ij}(\mathbf d)}-\sum_{j}\widehat P_D(j\mid i)\,a_{ij}(\mathbf d)\Big]\\
\text{s.t.}\ &\LM+\mathbf g_M^{\top}\mathbf d\le\varepsilon_M,\qquad \Lf+\mathbf g_f^{\top}\mathbf d\le\varepsilon_f,\\
&\mathbf s\,\mathrm e^{\mathbf d}\in\Sset,\qquad \|\mathbf d\|_\infty\le\Delta,
\end{aligned}
\label{eq:local-step}
\end{equation}
with $\Delta$ the trust radius. The objective of \eqref{eq:local-step} is a log-sum-exp of affine functions of $\mathbf d$ minus a term linear in $\mathbf d$, and its constraints are linear, so \eqref{eq:local-step} is a convex program. Its solution $\mathbf d^\star$ gives the trial scales $\mathbf s'=\mathbf s\,\mathrm e^{\mathbf d^\star}$ and the warm start $\mathbf z^\star(\mathbf s)+\mathbf U\operatorname{diag}(\mathbf s)\,\mathbf d^\star$, from which Algorithm~\ref{alg:hld} reaches $\mathbf z^\star(\mathbf s')$ in a few iterations. A filter on the losses at $\mathbf s'$ accepts or rejects the step and updates $\Delta$ \cite{fletcher2002filter}. The validation loss would vary with another draw of the validation origins. The search stops when the predicted decrease of $\LD$ falls below this variation, estimated by a paired bootstrap over the validation origins \cite{efronTibshirani1993}, since a smaller improvement cannot be distinguished from sampling noise.

\section{NUMERICAL RESULTS}
\label{sec:results}

\subsection{Data}
\label{sec:data}

We evaluate the proposed bi-level estimation framework on commuting trips in Los Angeles County, using three public datasets introduced as follows.

\textit{Commuting flows.} For the origin--destination flows and the destination attributes, we use LODES 2019 \cite{lehdLodes}, which records the census tract of residence and of workplace of 3.69 million workers of Los Angeles County.  We use 2,475 of the county's 2,498 census tracts because the remaining 23 tracts carry no or negligible commuting. The destination choice set $\mathcal J_i$ of an origin $i$ consists of the workplace tracts that receive at least 5 of its commuters. LODES data also gives the jobs and the resident workers of each tract by 20 industry sectors, from which the job counts $\widehat A_j$ and the sector match $\widehat\kappa_{ij}$ are computed. 

\textit{Mode shares.} For the mode level, we use the American Community Survey (ACS) aggregated five-year estimates for 2015--2019 \cite{acs2019}, which report the split $\widehat s_{im}$ of its workers across modes $m$ for each tract of residence $i$, i.e., driving alone, carpooling, transit, walking and bicycling. Since LODES carries no mode information, we construct the mode-specific counts as $\widehat T_{ijm}=\widehat T_{ij}\,\widehat s_{im}$ by assigning each origin tract's mode shares to all of its destinations and renormalizing over the modes available for each pair, with walking and bicycling only available within 5 and 15~km.

\textit{Network and traffic counts.} The road network is OpenStreetMap \cite{osm}, with free-flow speed and capacity as link attributes, and the transit network is built from GTFS feeds \cite{gtfs}. Driving alone and carpooling choose among three generated routes per pair, while transit, walking and bicycling move at fixed speeds of 22, 4.8 and 15~km/h on their own connections. PeMS \cite{pems} supplies average weekday counts at 1,593 mainline freeway detectors, summed over the morning commute window of 6--10~am. The PeMS counts include trips of all purposes, whereas the model loads commuting trips only, so the two differ by a common factor. The flow loss of Section~\ref{sec:model} removes this factor by subtracting the mean log ratio, so the counts need no rescaling. Each detector is matched to the sequence of consecutive network links it observes, 3.4 links on average, and its modeled flow is the mean over those links.

\textit{Problem scale.} We randomly split the 2,475 tracts in a 10:1 ratio into 2,250 estimation origins and 225 validation origins. The estimation origins carry 1.71 million commuters over 143,418 origin--destination pairs and 1.23 million routes, which with the mode and link states makes 2.7 million variables in \eqref{eq:inner}, and the validation origins carry 184,630 commuters over 15,006 pairs and 129,006 routes.

\textit{Specification.} We set $\widehat\tau_n=1$ and leave the estimation of the nest parameters to future work. For every road link, we set the BPR parameters to $\widehat a^m=0.15$ and $\widehat b^m=4$. The reference scales are $\mathbf s^{\mathrm{ref}}=(1,1,1)$. At these scales with $\widehat\tau_n=1$, the three logits collapse into one multinomial logit over destination, mode and route bundles with utility $V^D_{ij}+V^M_{ijm}+V^R_{ijmr}$, and \eqref{eq:objective} becomes the entropy of the joint pattern less the Beckmann potential, the classical combined model \cite{evans1976} with logit route choice, so the reference is the model before any scale is calibrated. Following the standard $\varepsilon$-constraint method \cite{haimes1971}, we set the tolerances in \eqref{eq:bilevel} to the losses at the reference scales $\mathbf{s}^{\mathrm{ref}}=(1,1,1)$: $\varepsilon_M=4.448$ and $\varepsilon_f=1.785$.

\subsection{Computational performance}
\label{subsec:performance}

All computations run on a server with two 96-core AMD EPYC 9655 processors and 1.5~TB of memory. As a computational baseline, we formulate the inner problem \eqref{eq:inner} as an exponential-cone program and solve it with the interior-point solver of MOSEK \cite{mosek} running on 184 parallel threads. The baseline runs the same outer search, except with cold starts in place of warm starts and finite differences in place of implicit differentiation, since MOSEK offers neither. Table~\ref{tab:timing} reports the total time of the calibration from $\mathbf s^{\mathrm{ref}}$ and the time and memory of one cold inner solve. HLD stops each inner solve by its own residual test, at which its destination loss is within $10^{-4}$ of that at the optimal solution of the baseline. On the full county, a cold inner solve takes HLD 7.5 minutes and 6.6~GB, against 194 minutes and 446~GB for the baseline. The memory of the baseline grows with its number of threads, so part of this memory comes from the 184 threads. With 16 threads, the baseline still needs 277~GB. HLD completes the calibration at every scale, in 0.74~h on the full county, whereas the baseline completes it only at 100 origins, in 8.1~h, and does not finish within 24~h from 200 origins on. The gap comes from the design of HLD, which warm-starts every inner solve after the first and differentiates its optimality conditions cheaply. The baseline cold-starts every solve and, with optimality conditions too large to differentiate, needs additional cold solves per outer step for the sensitivities.

\begin{table}[!t]
    \centering
    \footnotesize
    \setlength{\tabcolsep}{3pt}
    \renewcommand{\arraystretch}{1.05}
    \caption{Cost of calibration from $\mathbf s^{\mathrm{ref}}$}
    \label{tab:timing}
    \begin{threeparttable}
    \begin{tabular}{@{}rrrrrrrr@{}}
        \toprule
        \multicolumn{2}{c}{\multirow{2}{*}{Problem size}} & \multicolumn{2}{c}{\multirow{2}{*}{Total time (h)}} & \multicolumn{4}{c}{Per inner solve} \\
        \cmidrule(l){5-8}
        & & & & \multicolumn{2}{c}{Time (min)} & \multicolumn{2}{c}{Memory (GB)} \\
        \cmidrule(r){1-2}\cmidrule(lr){3-4}\cmidrule(lr){5-6}\cmidrule(l){7-8}
        Origins & Routes & HLD & Conic & HLD & Conic & HLD & Conic \\
        \midrule
        100 & 99,885 & 0.02 & 8.1 & 0.07 & 35 & 1.5 & 154 \\
        200 & 189,213 & 0.03 & -- & 0.26 & 61 & 1.9 & 159 \\
        400 & 348,180 & 0.15 & -- & 0.78 & 108 & 2.7 & 251 \\
        2,250 & 1,233,892 & 0.74 & -- & 7.5 & 194 & 6.6 & 446 \\
        \bottomrule
    \end{tabular}
    \begin{tablenotes}[flushleft]\scriptsize
    \item[] -- marks a calibration not finished within the 24~h cap. Time and memory per inner solve are those of one cold solve at $\mathbf s^{\mathrm{ref}}$.
    \end{tablenotes}
    \end{threeparttable}
\end{table}

\begin{figure}[h]
    \centering
    \includegraphics[width=.74\columnwidth]{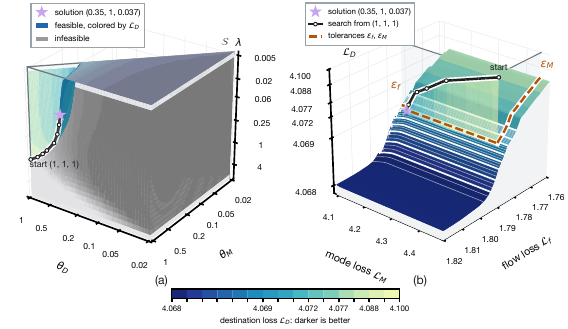}
    \caption{The outer search on the county instance: (a) the trajectory of the scales in $\Sset$, (b) $\LD$ against the constrained losses $\LM$ and $\Lf$ along it. The $\LD$ scale is stretched near its minimum.}
    \label{fig:outer3d}
\end{figure}

\subsection{The outer search on the county instance}
\label{subsec:scales}

Fig.~\ref{fig:outer3d} shows the outer search process from the reference scales $\mathbf s^{\mathrm{ref}}=(1,1,1)$, where both tolerances hold with equality and the objective is $\LD=4.097$. It ends at $\mathbf s^\star=(0.35,1,0.037)$ with the network tolerance active and $\LD=4.069$. The outer search thus extends \cite{kim2026estimate} by calibrating the scales that reconcile the destination, mode and network observations within one model.

\subsection{Parameter estimates and predictive performance}
\label{subsec:estimates}

Table~\ref{tab:parameters} reports the solution of the inner problem \eqref{eq:inner}, the recovered taste coefficients of the destination and mode levels, at the calibrated scales $\mathbf s^\star=(0.35,1,0.037)$ alongside the solution at the reference scales $\mathbf s^{\mathrm{ref}}=(1,1,1)$. Since the taste coefficients themselves depend on the scale, we report their products with the corresponding scales, which determine the choice probabilities and are directly comparable across scale specifications. The recovered coefficients in our calibrated model have the expected signs and plausible magnitudes: doubling the jobs at a workplace raises its odds of being chosen by a factor of 1.6, a workplace 20~km away has 40\% lower odds than one 10~km away, and a 0.1 increase in the sector match multiplies the odds by 1.6. At the mode level, an extra hour of travel time lowers a mode's odds by 26\%. Most importantly, at the reference scales ($\lambda=1$), the time coefficient is positive, showing that incorrectly fixed scales can yield behaviorally implausible estimates and motivating the outer calibration.

\begin{table}[h!]
    \centering\scriptsize\setlength{\tabcolsep}{3pt}\renewcommand{\arraystretch}{0.95}
    \caption{Estimates and predictive performance at $\mathbf s^{\mathrm{ref}}$ and  $\mathbf s^\star$.}
    \label{tab:parameters}
    \begin{tabular}{@{}lrr@{}}
        \toprule
        $(\theta_D,\theta_M,\lambda)$ & Reference $\mathbf s^{\mathrm{ref}}$ & Calibrated $\mathbf s^\star$ \\
        \midrule
        \emph{Estimates, destination level} & $\theta_D\boldsymbol\beta^D$ & $\theta_D\boldsymbol\beta^D$ \\
        distance $\widehat d_{ij}$ (km) & $-0.015$ & $-0.023$ \\
        log jobs $\ln(1+\widehat A_j)$ & $0.688$ & $0.641$ \\
        log distance $\ln\widehat d_{ij}$ & $-0.465$ & $-0.413$ \\
        sector match $\widehat\kappa_{ij}$ & $5.036$ & $4.438$ \\
        \emph{Estimates, mode level} & $\theta_M\boldsymbol\beta^M$ & $\theta_M\boldsymbol\beta^M$ \\
        carpool constant & $-2.100$ & $-2.100$ \\
        walk constant & $-2.501$ & $-2.508$ \\
        bicycle constant & $-3.792$ & $-3.776$ \\
        transit constant & $-1.893$ & $-1.863$ \\
        travel time $\widehat t_{ijm}$ (h) & $0.618$ & $-0.305$ \\
        \midrule
        \multicolumn{3}{@{}l}{\emph{Predictive performance}} \\
        $\LD$ & $4.097$ & $\mathbf{4.069}$ \\
        $\LM$ & $4.448$ & $\mathbf{4.118}$ \\
        $\Lf$ & $1.785$ & $1.785$ \\
        $R^2$ & $0.609$ & $\mathbf{0.624}$ \\
        CPC & $0.785$ & $\mathbf{0.804}$ \\
        \bottomrule
    \end{tabular}
\end{table}

\begin{figure}[!t]
    \centering
    \includegraphics[width=.66\columnwidth]{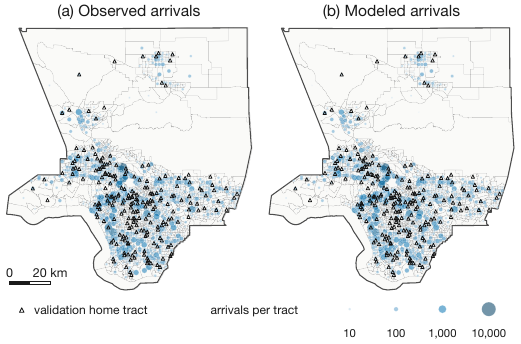}
    \caption{Commuting destinations by workplace tract from the validation origins, (a) observed and (b) modeled at the calibrated scales $\mathbf s^\star$}
    \label{fig:maps}
\end{figure}

Table~\ref{tab:parameters} also reports the three validation losses, $\mathcal L_D$, $\mathcal L_M$, and $\mathcal L_f$, together with two diagnostics of the destination prediction on the validation origins: the $R^2$ of the modeled against the observed pair flows, and the common part of commuters (CPC) \cite{lenormandEtAl2016}, with $\mathrm{CPC}=\frac{\sum_{i\in\mathcal I_{\mathrm{val}}}\sum_j\min(T_{ij},\widehat T_{ij})}{\sum_{i\in\mathcal I_{\mathrm{val}}}\sum_j\widehat T_{ij}},$
the share of the validation commuters that the model places on the same pairs as the data. The calibrated model is applied to the validation tracts with the estimated coefficients. Its destination loss $\LD$ falls from 4.097 at the reference scales to 4.069, its mode loss $\LM$ falls from 4.448 to 4.118, and its flow loss $\Lf$ stays at 1.785. Both validation diagnostics improve, $R^2$ from 0.609 to 0.624 and CPC from 0.785 to 0.804. Fig.~\ref{fig:maps} shows that the modeled arrivals closely reproduce the observed spatial concentration, with only slightly more dispersion around the largest centers.

\begin{figure}[h]
    \centering
    \includegraphics[width=.72\columnwidth]{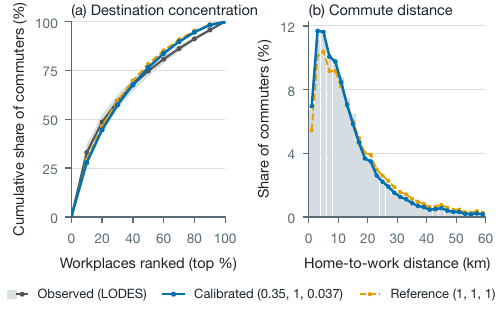}
    \caption{Destination patterns (a) cumulative share of a tract's commuters over its workplaces ranked by visits (b) distribution of home-to-work distance.}
    \label{fig:distributions}
\end{figure}

Fig.~\ref{fig:distributions} shows two further patterns on the validation origins. In panel (a), the cumulative share of a tract's commuters over its workplaces follows the observed mean curve to within 5.3 percentage points. In panel (b), although our calibrated model is fitted only to the mean distance and the mean log distance of the estimation origins, it reproduces the full distribution of commute distances well. The calibrated model follows the observed distribution over its whole range, including the peak below 10~km, whereas the reference scales place too few commuters on trips under 10~km and too many beyond 25~km, which lengthens their mean commute from the observed 13.3~km to 15.2~km.

\section{CONCLUSION}

In this paper, we propose a bi-level program that estimates a combined destination, mode, and route choice model with congestion, from aggregate data. It recovers the taste coefficients, the equilibrium travel pattern, and the scale parameters jointly, and we test it at the scale of Los Angeles County with accurate prediction on the validation origins. There are still interesting future directions. First, the inputs of our model are public aggregate data, and whether the behavioral parameters recovered from them agree with those conventionally estimated from individual-level household travel surveys is an open question. Second, activity-based models characterize current travel behavior in rich detail through individual activity patterns and traveler heterogeneity, whereas our trip-based parametric model is parsimonious, interpretable and can be carried into future scenarios. A hybrid framework, in which an activity-based model characterizes the present system and the parametric combined model estimates behavior and predicts future scenarios, is a natural next step.

\section*{ACKNOWLEDGMENT}

This study was made possible with funding received by the University of California Institute of Transportation Studies through the Road Repair and Accountability Act of 2017 (Senate Bill~1). This paper is disseminated under the sponsorship of the State of California in the interest of information exchange and does not necessarily reflect the official views or policies of the State of California.

\bibliographystyle{IEEEtran}
\bibliography{references}

\end{document}

%% file: inputs/bilevel_framework.tex
\begin{figure}[!t]
    \centering
    \begin{tikzpicture}[
        x=1cm, y=1cm,
        alt/.style   = {draw, ellipse, minimum width=13mm, minimum height=7mm,
                        inner sep=0pt, font=\small, line width=0.5pt},
        band/.style  = {draw, dashed, rounded corners=5pt, black!45, line width=0.5pt},
        lvl/.style   = {font=\small, align=left, anchor=west},
        pr/.style    = {font=\small, inner sep=2pt, fill=white},
        lk/.style    = {line width=0.5pt, black!75},
        up/.style    = {line width=0.6pt, black!55, -{Stealth[length=2.2mm]}},
        uplab/.style = {font=\footnotesize, rotate=90, inner sep=2pt},
        box/.style   = {draw, rounded corners=6pt, line width=0.6pt, inner sep=4pt},
        stp/.style   = {draw, rounded corners=3pt, line width=0.5pt, font=\footnotesize,
                        align=center, inner sep=2.5pt, fill=white, text width=58mm},
        flow/.style  = {line width=0.7pt, -{Stealth[length=2.4mm]}},
        flab/.style  = {font=\scriptsize, align=center, fill=white, inner sep=1.5pt},
        ttl/.style   = {font=\footnotesize\bfseries, anchor=north west, inner sep=3pt}]

      \begin{scope}[scale=0.5, every node/.append style={transform shape}]
        \draw[band] (-3.3,0.75)  rectangle (7.1,-2.45);
        \draw[band] (-3.3,-2.95) rectangle (7.1,-4.45);
        \draw[band] (-3.3,-5.00) rectangle (7.1,-6.70);
        \node[alt] (O)  at (0,0)       {$\widehat O_i$};
        \node[alt] (D1) at (-2.3,-1.7) {};
        \node[alt] (Dj) at (0,-1.7)    {$T_{ij}$};
        \node[alt] (D3) at (2.3,-1.7)  {};
        \node[alt] (M1) at (-2.3,-3.7) {};
        \node[alt] (Mm) at (0,-3.7)    {$T_{ijm}$};
        \node[alt] (M3) at (2.3,-3.7)  {};
        \node[alt] (R1) at (-2.3,-5.85){};
        \node[alt] (Rr) at (0,-5.85)   {$T_{ijmr}$};
        \node[alt] (R3) at (2.3,-5.85) {};
        \draw[lk] (O)--(D1);  \draw[lk] (O)--(Dj);  \draw[lk] (O)--(D3);
        \draw[lk] (Dj)--(M1); \draw[lk] (Dj)--(Mm); \draw[lk] (Dj)--(M3);
        \draw[lk] (Mm)--(R1); \draw[lk] (Mm)--(Rr); \draw[lk] (Mm)--(R3);
        \node[pr] at (0,-0.85) {$P_D(j\mid i)$};
        \node[pr] at (0,-2.70) {$P_M(m\mid i,j)$};
        \node[pr] at (0,-4.72) {$P_R(r\mid i,j,m)$};
        \draw[up] (-3.9,-5.85) -- (-3.9,-3.70);
        \node[uplab] at (-4.22,-4.78) {$S^R_{ijm}$};
        \draw[up] (-3.9,-3.70) -- (-3.9,-0.85);
        \node[uplab] at (-4.22,-2.28) {$S^M_{ij}$};
        \node[lvl] at (3.4,-0.85) {Destination choice\\(Multinomial logit)};
        \node[lvl] at (3.4,-3.70) {Mode choice\\(Nested logit)};
        \node[lvl] at (3.4,-5.85) {Route choice\\(Path-size logit)};
        \coordinate (hNW) at (-4.7,1.6);
        \coordinate (hSE) at (7.3,-6.9);
      \end{scope}
      \node[box, fit=(hNW)(hSE)] (inner) {};
      \node[ttl] at (inner.north west) {Inner problem \eqref{eq:inner}: estimation at fixed scale $\mathbf s$};

      \node[stp, above=3.5mm of inner] (cand)
        {Candidate scales $\mathbf s=(\theta_D,\theta_M,\lambda)\in\Sset$};
      \node[stp, below=3.5mm of inner] (pred)
        {Evaluation on $\mathcal I_{\mathrm{val}}$:\ losses $\LD$, $\LM$, $\Lf$};
      \node[stp, below=2.5mm of pred] (sel)
        {Selection:\ $\operatorname*{arg\,min}_{\mathbf s}\LD$\ \\s.t.\ $\LM\le\varepsilon_M$,\ $\Lf\le\varepsilon_f$};
      \draw[flow] (cand) -- (inner) node[flab, midway, right=1mm] {$\mathbf s$};
      \draw[flow] (inner) -- (pred) node[flab, midway, right=1mm] {$\boldsymbol\beta^\star(\mathbf s)$, $(\mathbf p^\star,\mathbf f^\star)$};
      \draw[flow] (pred) -- (sel);
      \draw[flow] (sel.west) -- ++(-4mm,0) |- (cand.west);
      \coordinate (top)   at ([yshift=5.5mm]cand.north);
      \coordinate (loopL) at ([xshift=-7.5mm]sel.west);
      \node[box, line width=0.9pt, fit=(cand)(inner)(pred)(sel)(top)(loopL), inner sep=3pt] (outer) {};
      \node[ttl] at (outer.north west) {Outer problem \eqref{eq:bilevel}: scale selection};
    \end{tikzpicture}
        \caption{Bi-level structure of the estimation problem.}
    \label{fig:choice-hierarchy}
\end{figure}